# Dimensionality Reduction in Stochastic Complex Dynamical Networks

**Authors:** Chengyi Tu[1,2*], Jianhong Luo[1], Xuwei Pan[1*]

[1]School of Economics and Management, Zhejiang Sci-Tech University; Hangzhou, 310018, China.

[2]Department of Environmental Science, Policy, and Management, University of California, Berkeley; Berkeley, 94720, USA.

*Corresponding author. Email: chengyitu1986@gmail.com and panxw@zstu.edu.cn

# Abstract

Complex systems often exhibit diverse dynamical behaviors in high-dimensional spaces that depend on various factors. Dimensionality reduction is a powerful tool for analyzing and understanding complex systems, aiming to find a low-dimensional representation of the complex system that preserves its essential features and reveals its underlying mechanisms and long-term dynamics. However, most existing methods for dimensionality reduction are limited to deterministic systems and cannot account for the stochastic effects that are ubiquitous in real-world complex networks. Here we develop a general analytical framework for dimensionality reduction of stochastic complex dynamical networks that can capture the essential features and long-term dynamics of the original system in a low-dimensional effective equation. The effective equation is a function of a set of effective parameters that are associated with specific system states and determine the network's dynamical behavior. We show that the standard deviation of the effective equation can be used to analyze the dynamic behavior and possible convergence of the stochastic complex dynamical network. Our framework can be applied to various types of stochastic complex dynamical networks and can reveal the underlying mechanisms and emergent phenomena of these systems.

**Keywords**: dimensionality reduction; stochastic complex dynamical networks; deterministic and stochastic effect

# Introduction

Dimensionality reduction of complex systems is a topic of growing interest across a wide range of scientific disciplines, such as biology, data science, chemistry, microelectronics and socioeconomics [1-6]. Complex systems often exhibit rich and diverse dynamical behaviors that depend on various factors, such as interaction structures, parameter configurations and initial conditions [7, 8]. These behaviors are embedded in high-dimensional spaces that pose significant challenges for analysis and understanding. Dimensionality reduction aims to find a low-dimensional representation of the complex system that preserves its essential features and reveals its underlying mechanisms and long-term dynamics. Such a representation can facilitate the control and design of the system, as well as the prevention of undesirable transitions. Dimensionality reduction is therefore a valuable tool for both practical applications and theoretical investigations.

Despite the considerable efforts devoted to dimensionality reduction of complex systems [9], many challenges remain, especially for complex dynamical networks. Previous methods, such as the one-dimensional approach [10] and the critical slowing down theory [11-14], have limited applicability or generality. Gao et al. made the first attempt to apply dimensionality reduction to complex dynamical networks by developing an analytical tool that can extract the effective parameters of a high-dimensional networked system using mean-field approximations [15]. However, Tu et al. found a new condition that severely restricts the validity of Gao's framework [16]. They also extended Gao's framework to account for the heterogeneity in node-specific self- and coupling-dynamics, and applied it to study the collapse or functioning of various networked systems [17]. Moreover, to deal with arbitrary network structures, Laurence et al. proposed a polynomial approximation framework to reduce the complexity of dynamical networks [18, 19]. Wu et al. introduced an improved dimension reduction method based on information entropy to predict network resilience [20].

However, these methods are limited to complex dynamical networks without stochastic effects. In reality, the dynamic behavior and signal transmission are often subject to various random noises that may arise from random fluctuations in the release of uncertainties or environmental perturbations. Therefore, stochastic models are essential for capturing more realistic dynamical behaviors of complex networks, and many studies have explored the complex dynamical networks with stochastic effects. To overcome the theoretical and practical challenges of dimensionality reduction for stochastic complex dynamical networks, a general framework that can account for the stochastic effects is urgently needed.

In this work, we present a general analytical framework for dimensionality reduction of stochastic complex dynamical networks. The framework can reduce the order parameter space of the network to a function of a set of effective parameters, which are defined as those that determine the network's dynamical behavior and are associated with specific system states. The resulting effective equation is a low-dimensional representation of the original high-dimensional stochastic complex dynamical network that preserves its essential dynamical features. The framework also enables us to analyze the dynamic behavior and possible convergence of the stochastic complex dynamical network by using the standard deviation of the effective equation. Our framework can be applied to various types of stochastic complex dynamical networks, such as biology [21], ecology [22], epidemiology [23], neuroscience [24] and so on, and can reveal the underlying mechanisms and emergent phenomena of these systems.

## Methods

We consider a stochastic complex dynamical network of $N$ nodes whose states $\mathbf{x}=(x_1,\ldots,x_N)^T$ evolve according to the following dynamic equation:

$$dx_i = \left( F_i(x_i) + \sum_{j}^{N} A_{ij} G_i(x_i, x_j) \right) dt + H_i(x_i) dW_i \quad (1)$$

where $F_i(x_i)$ denotes the self-dynamics of node $i$, $G_i(x_i, x_j)$ denotes the coupling dynamics between node $i$ and its neighbor $j$, $A_{ij}$ denotes the interaction strength between node $i$ and node $j$, and $H_i(x_i) W_i$ denotes the diffusion term of the stochastic differential equation for node $i$ with $W_i$ being a Wiener process.

From Eq. (1), we can see that the dynamics of each node are influenced by three factors: the self-dynamics $F_i(x_i)$, the interaction with its nearest neighbors (captured by the interaction network $\mathbf{A}$ and the coupling dynamics $G_i(x_i, x_j)$) and the stochastic term $H_i(x_i) W_i$. To analyze the dynamical behavior of the stochastic complex dynamical network, we introduce a mean-field operator $\mathcal{L}(\mathbf{x}) = \frac{1}{N}\sum_{j=1}^{N} s_j^{out} x_j \bigg/ \frac{1}{N}\sum_{j=1}^{N} s_j^{out} = \frac{\langle \mathbf{s}^{out} \cdot \mathbf{x} \rangle}{\langle \mathbf{s}^{out} \rangle}$, where $\mathbf{s}^{out} = (s_1^{out},\ldots,s_N^{out})^T$ is the vector of the out degree of the matrix $\mathbf{A}$ [15], and use the weighted average node state as an effective state of the stochastic complex dynamical network. This operator $\mathcal{L}$ is feasible for both linear time-invariance (LTI) functions and the Hadamard product, i.e., $\mathcal{L}(a\mathbf{x}+b\mathbf{y}) = a\mathcal{L}(\mathbf{x}) + b\mathcal{L}(\mathbf{y})$ and $\mathcal{L}(\mathbf{x} \circ \mathbf{y}) \approx \mathcal{L}(\mathbf{x})\mathcal{L}(\mathbf{y})$ [17].

Under some assumptions, we can simplify Eq. (1) using the mean-field operator $\mathcal{L}$. First, if the network $\mathbf{A}$ has low degree correlation (meaning that the nodes have similar neighborhoods), then we can approximate the interaction term as $\sum_{j}^{N} A_{i,j} G_i(x_i, x_j) \approx s_i^{in} \mathcal{L}(G_i(x_i, \mathbf{x}))$. Second, if $G_i(x_i, x_j)$ is linear in $x_j$ or the vector $\mathbf{x}$ has small standard deviation, then we can approximate the mean-field term as $\mathcal{L}(G_i(x_i, \mathbf{x})) \approx G_i(x_i, \mathcal{L}(\mathbf{x})) = G_i(x_i, x_{eff})$, where $x_{eff} = \mathcal{L}(\mathbf{x})$ is the effective state of the network. Therefore, we can rewrite Eq. (1) as $dx_i \approx \left( F_i(x_i) + s_i^{in} G_i(x_i, x_{eff}) \right) dt + H_i(x_i) dW_i$, and in vector form as $d\mathbf{x} = \left( \mathbf{F}(\mathbf{x}) + \mathbf{s}^{in} \circ \mathbf{G}(\mathbf{x}, x_{eff}) \right) dt + \mathbf{H}(\mathbf{x}) \circ d\mathbf{W}$, where $\mathbf{F}(\mathbf{x}) = \left( F_1(x_1),\ldots,F_N(x_N) \right)^T$ ,

$$\mathbf{G}(\mathbf{x}, x_{e\!f\!f}) = \big(G_1(x_1, x_{e\!f\!f}), \ldots, G_N(x_N, x_{e\!f\!f})\big)^T, \quad \mathbf{s}^{in} = (s_1^{in}, \ldots, s_N^{in})^T, \quad \mathbf{H}(\mathbf{x}) = \big(H_1(x_1), \ldots, H_N(x_N)\big)^T \text{ and}$$

$$\mathbf{W} = (W_1, \ldots, W_N)^T.$$

We can further simplify the self-dynamics term $\mathbf{F}(\mathbf{x})$ if we assume that each $F_i(x_i)$ is a linear combination of $m$ subfunctions, i.e., $F_i(x_i) = b_{i,1} f_1(x_i) + b_{i,2} f_2(x_i) + \cdots + b_{i,m} f_m(x_i)$. Using the linearity property of the operator $\mathcal{L}$, we can write $\mathcal{L}(\mathbf{F}(\mathbf{x})) = \mathcal{L}\begin{pmatrix} F_1(x_1) \\ \vdots \\ F_N(x_N) \end{pmatrix} = \mathcal{L}\begin{pmatrix} b_{1,1} f_1(x_1) \\ \vdots \\ b_{N,1} f_1(x_N) \end{pmatrix} + \ldots + \mathcal{L}\begin{pmatrix} b_{1,m} f_m(x_1) \\ \vdots \\ b_{N,m} f_m(x_N) \end{pmatrix}$. Then, using the Hadamard product property of the operator $\mathcal{L}$, we can approximate $\begin{aligned}\mathcal{L}(\mathbf{F}(\mathbf{x})) &= \mathcal{L}(B^1 \circ f_1(\mathbf{x})) + \cdots + \mathcal{L}(B^m \circ f_m(\mathbf{x})) \\ &\approx \mathcal{L}(B^1)\mathcal{L}(f_1(\mathbf{x})) + \cdots + \mathcal{L}(B^m)\mathcal{L}(f_m(\mathbf{x}))\end{aligned}$ where $B^k = (b_{1,k}, \ldots, b_{N,k})^T$ is the $k$-th column of matrix $\mathbf{B}$. Since the subfunctions are the same for all nodes, we can further approximate $\begin{aligned}\mathcal{L}(\mathbf{F}(\mathbf{x})) &\approx \mathcal{L}(B^1) f_1(\mathcal{L}(\mathbf{x})) + \ldots + \mathcal{L}(B^m) f_m(\mathcal{L}(\mathbf{x})) \approx \\ \mathcal{L}(B^1) f_1(x_{e\!f\!f}) &+ \ldots + \mathcal{L}(B^m) f_m(x_{e\!f\!f}) = \sum_{k=1}^m B_{e\!f\!f}^k f_k(x_{e\!f\!f})\end{aligned}$ where $B_{e\!f\!f}^k = \mathcal{L}(B^k)$. We can apply a similar procedure to simplify the coupling dynamics term $\mathbf{G}_i(\mathbf{x}_i, \mathbf{x}_j)$ and the stochastic term $H_i(x_i) W_i$. If we assume that each $G_i(x_i, x_j)$ is a linear combination of $n$ subfunctions, i.e., $G_i(x_i, x_j) = c_{i,1} g_1(x_i, x_j) + \cdots + c_{i,n} g_n(x_i, x_j)$, then we can use the linearity and Hadamard product properties of the operator $\mathcal{L}$ to approximate

$$\mathcal{L}(\mathbf{G}(\mathbf{x}, x_{e\!f\!f})) \approx \mathcal{L}(C^1) g_1(x_{e\!f\!f}, x_{e\!f\!f}) + \ldots + \mathcal{L}(C^n) g_n(x_{e\!f\!f}, x_{e\!f\!f}) = \sum_{l=1}^n C_{e\!f\!f}^l g_l(x_{e\!f\!f}, x_{e\!f\!f}) \quad \text{where} \quad C_{e\!f\!f}^l = \mathcal{L}(C^l) \quad \text{and}$$

$C^l = (c_{1,l}, \ldots, c_{N,l})^T$ is the $l$-th column of matrix $\mathbf{C}$. Similarly, if we assume that each $H_i(x_i)$ is a linear combination of $t$ subfunctions, then we can approximate $\mathcal{L}(\mathbf{H}(\mathbf{x})) \approx \sum_{s=1}^t D_{e\!f\!f}^s h_s(x_{e\!f\!f})$ where $D_{e\!f\!f}^s = \mathcal{L}(D^s)$ and $D^s = (d_{1,s}, \ldots, d_{N,s})^T$ is the $s$-th column of matrix $\mathbf{D}$. Additionally, we can approximate $\mathcal{L}(\mathbf{W}) \approx W_{e\!f\!f}$ where $W_{e\!f\!f}$ is a Wiener process. Finally, we can obtain $\mathcal{L}(\mathbf{H}(\mathbf{x}) \circ \mathbf{W}) = \sum_{s=1}^t D_{e\!f\!f}^s h_s(x_{e\!f\!f}) W_{e\!f\!f}$. Applying the operator $\mathcal{L}$ to both sides of the vector equation, we get the effective equation

$$\begin{aligned}dx_{e\!f\!f} &= \mathcal{L}\big(\big(\mathbf{F}(\mathbf{x}) + \mathbf{s}^{in} \circ \mathbf{G}(\mathbf{x}, \mathcal{L}(\mathbf{x}))\big) dt + \mathbf{H}(\mathbf{x}) \circ d\mathbf{W}\big) \\ &\approx \left(\sum_{k=1}^m B_{e\!f\!f}^k f_k(x_{e\!f\!f}) + A_{e\!f\!f} \sum_{l=1}^n C_{e\!f\!f}^l g_l(x_{e\!f\!f}, x_{e\!f\!f})\right) dt + \sum_{s=1}^t D_{e\!f\!f}^s h_s(x_{e\!f\!f}) dW_{e\!f\!f}\end{aligned}.$$

In some cases, the self-dynamics function $F_i(x_i)$ may not be a linear combination of $m$ subfunctions or may differ for each node. In such cases, we can use Chebyshev polynomials to approximate $F_i(x_i)$ as $\sum_{k=1}^{m} b_{i,k} x^{(k-1)}$, minimizing the approximation error [25, 26]. Then, we can apply the operator $\mathcal{L}$ to get $\mathcal{L}(\mathbf{F}(\mathbf{x}[t])) = \sum_{k=1}^{m} B_{eff}^{k} x_{eff}^{(k-1)}$. Similarly, we can use Chebyshev polynomials to approximate $G_i(x_i, x_j)$ as $\sum_{p,q=1}^{n/2} d_{p,q} x_i^{(p-1)} x_j^{(q-1)}$ and apply the operator $\mathcal{L}$ to get $\mathcal{L}(\mathbf{G}(\mathbf{x}, x_{eff})) = \sum_{l=1}^{n} C_{eff}^{l} x_{eff}^{(l-1)}$ where $c_{i,l}$ collects all terms $d_{p,q}$ such that $l = p + q - 1$. We can do the same for $H_i(x_i)$ and get $\mathcal{L}(\mathbf{H}(\mathbf{x})) = \sum_{s=1}^{t} D_{eff}^{s} x_{eff}^{(s-1)}$. Finally, we obtain the effective equation of the stochastic complex dynamical network [17]

$$dx_{eff} \approx \left( \sum_{k=1}^{m} B_{eff}^{k} x_{eff}^{(k-1)} + A_{eff} \sum_{l=1}^{n} C_{eff}^{l} x_{eff}^{(l-1)} \right) dt + \sum_{s=1}^{t} D_{eff}^{s} x_{eff}^{(s-1)} dW_{eff} \quad (2)$$

The effective equation, Eq. (2), is a low-dimensional approximation of the stochastic complex dynamical network, Eq. (1), so we can use its dynamical behavior to infer the collective dynamical behavior of the stochastic complex dynamical network. If the stochastic complex dynamical network has small fluctuations, then the effective equation also has small fluctuations, and the standard deviation of its realizations is small. This means that the diffusion term is negligible compared to the drift term and the stochastic effect is not significant. Moreover, if the stochastic complex dynamical network converges, then the trajectory of the standard deviation of the effective equation decreases after an initial jump. Therefore, the standard deviation of the effective equation is a useful indicator for the dynamic feature of the stochastic complex dynamical network: A small value indicates that the stochastic effect is not important for the stochastic complex dynamical network and a decreasing trend indicates convergence of the stochastic complex dynamical network.

# Results

## Dynamic feature of multivariate Ornstein-Uhlenbeck process

One of the most common types of stochastic processes that models the random motion of a particle under the influence of friction is the Ornstein-Uhlenbeck process [27]. It has many applications in finance, physics, and biology. We use a multivariate Ornstein-Uhlenbeck process as an example to test the validity of our framework. The process is described by the equation:

$$d\mathbf{x} = -\mathbf{A}\mathbf{x}dt + \mathbf{B}d\mathbf{W} \quad (3)$$

where **A** is a matrix and **B** is a vector with constant elements $B_i = \epsilon$ representing the stochastic strength and **W** is a vector of Wiener processes. We want to answer some questions about this process: How does the drift term affect the general trend? How weak should the diffusion term be for convergence to occur? How small should the stochastic strength be? We apply our framework to address these questions. The effective equation of Eq. (3) is

$$dx_{eff} = -A_{eff} x_{eff} dt + B_{eff} dW_{eff} \quad (4)$$

where $x_{eff} = \mathcal{L}(\mathbf{x})$, $A_{eff} = \mathcal{L}(\mathbf{s}^{in})$, $B_{eff} = \mathcal{L}(\mathbf{B}) = \epsilon$ and $W_{eff} = \mathcal{L}(\mathbf{W})$.

To validate our framework, we conduct numerical simulations to examine how the stochastic strength affects the dynamical behavior of the stochastic complex dynamical network. We set the parameters as follows: $N = 50$, $A_{ij}$ are drawn from a distribution with mean $\mu_A = 0.001$ and standard deviation $\sigma_A = \mu_A / 3$ and stochastic strength $\epsilon$ varies from 10^-1.5 to 10^0.5. We initialize $\mathbf{x}(t = 0)$ with random values between 0 and 1. The results are shown in Fig. 1. We observe that when the stochastic strength of the original stochastic dynamical network is small\large (Fig. 1 b\d), the standard deviation of the realizations of the effective equation is also small\large (Fig. 1 c\e), indicating that the diffusion\drift term dominates the dynamics. Therefore, the effective equation is a good approximation of the stochastic complex dynamical network, especially for determining whether the deterministic effect is stronger than the stochastic effect.

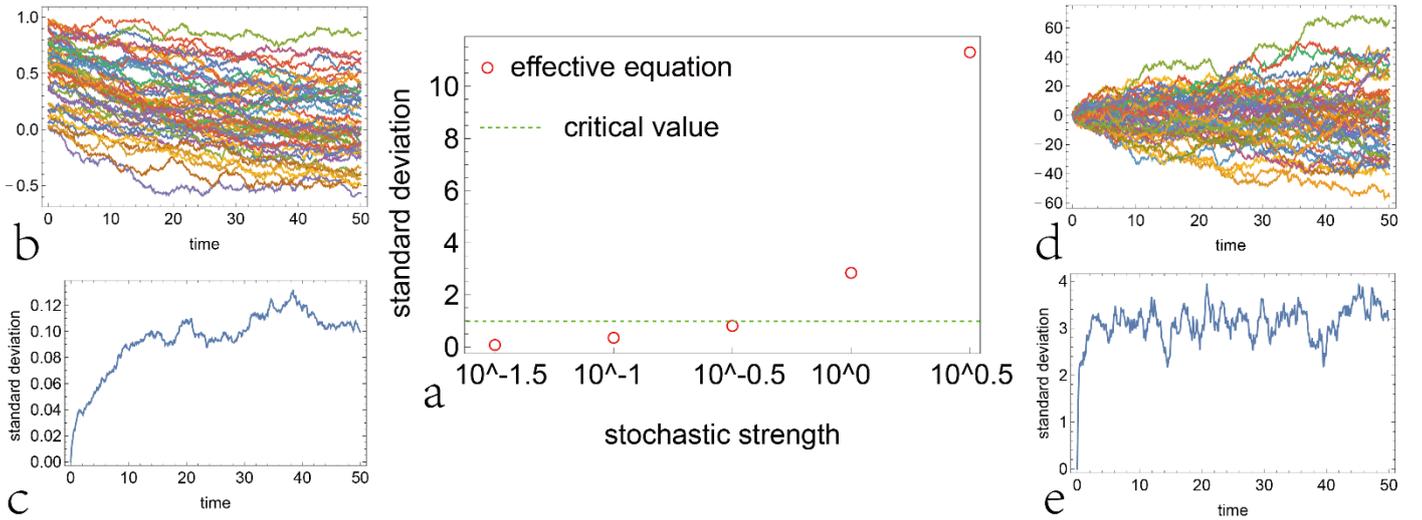

**Figure 1**: The effect of stochastic strength on the dynamics of the stochastic complex dynamical network and the effective equation. (a) The standard deviation of the realizations of the effective equation for different values of stochastic strength: 10^-1.5, 10^-1, 10^-0.5, 10^0 and 10^0.5. (b-c) A sample realization of the stochastic complex dynamical network with stochastic strength 10^-1.5 and the corresponding trajectory of the standard deviation of the effective equation. (d-e) A sample realization of the stochastic complex dynamical network with stochastic strength 10^0.5 and the corresponding trajectory of the standard deviation of the effective equation. We repeat each simulation 50 times.

# Dynamic feature of the GLV dynamics with empirical networks

In this section, we apply our framework to a generalized Lotka-Volterra (GLV) dynamics [28-30] with stochastic effect. The GLV dynamics are a set of first-order, nonlinear differential equations that are widely used in community ecology to model the dynamics of the populations of interacting species. The interactions between species can be represented by a network with different types of links, such as mutualism (++), commensalism (+0), competition (–), amensalism (-0), and predation (±). The GLV dynamics with stochastic effect for species $i$ are given by

$$dx_i = \left( \alpha_i x_i + \sum_{j=1}^{N} x_i A_{ij} x_j \right) dt + \epsilon x_i dW_i \quad (5)$$

where $N$ is the number of species in the community, $x_i$ is the population size of species $i$, $\alpha_i$ is the intrinsic growth rate of species $i$, $A_{ij}$ is the interaction strength between species $i$ and $j$, $W_i$ is a Wiener process for species $i$ and $\epsilon$ is a constant representing the stochastic strength. The effective equation of Eq. (5) is

$$dx_{eff} = \left( \alpha_{eff} x_{eff} + A_{eff} x_{eff}^2 \right) dt + \epsilon x_{eff} dW_{eff} \quad (6)$$

where $x_{eff} = \mathcal{L}(\mathbf{x})$, $A_{eff} = \mathcal{L}(\mathbf{s}^{in})$, $\alpha_{eff} = \mathcal{L}(\boldsymbol{\alpha})$ and $W_{eff} = \mathcal{L}(\mathbf{W})$.

As an example, we consider a community of $N$ species that consists of $N_p$ plants and $N_a$ animals such as insects that act as pollinators, with $N = N_p + N_a$ [31-35]. We use the vector $\mathbf{x} = \{x_1^p, x_2^p, \ldots, x_{N_p}^p, x_{N_p+1}^a, \ldots, x_N^a\}$ to represent the population size of each species in the community, where $x_i^p$ and $x_j^a$ denote the abundances of the $i$-th plant species and $j$-th animal species. The vector $\boldsymbol{\alpha}$ is the intrinsic growth rate of each species. The interaction matrix $\mathbf{A}$ has four blocks: two of them capture the direct competitive interactions between plants ($\Omega_{pp}$) and animals ($\Omega_{aa}$), while the other two describe the mutualistic interactions between plants and animals ($\Gamma_{pa}$) and vice versa ($\Gamma_{ap}$). The structure of the interaction matrix is given by $\begin{bmatrix} \Omega_{pp} & \Gamma_{pa} \\ \Gamma_{ap} & \Omega_{aa} \end{bmatrix}$, which has both positive and negative elements and a nonzero correlation [32]. We set the matrices $\Omega_{pp}$ and $\Omega_{aa}$ to represent the competition coefficients $\beta$ between plants and animals. We draw $\beta$ from a uniform distribution with a maximum of -0.001 and a mean value of $-1/S^p$ and $-1/S^a$. The matrices $\Gamma_{pa}$ and $\Gamma_{ap}$ capture the mutualistic interactions between plants and animals. We use a trade-off function that defines the mutualistic dependence between species $j$ and $i$ as a function of their degree: $\gamma_{ij} = \gamma y_{ij} / k_i$, where $\gamma$ is drawn from a normal distribution with mean $\mu_\gamma$ and standard deviation $\sigma_\gamma = |\mu_\gamma/3|$, $k_i$ is the degree of species $i$ and $y_{ij} = 1$ if species $i$ and $j$ interact; otherwise, it is zero [32]. We construct the interaction matrix using the Web-of-Life database (http://www.web-of-life.es), which contains empirical plant-pollinator networks from different regions of the world. We select 134 networks with no more than 200 species to reduce the computational cost of our analysis [3]. We set the parameters as follows: $\mu_\alpha = 1, \sigma_\alpha = |\mu_\alpha/3|$ and $\mu_\gamma = 0.4, \sigma_X = |\mu_\gamma/3|$. We initialize $\mathbf{x}(t=0)$ with random values between 0 and 1.

We use numerical simulations to investigate how the stochastic strength affects the dynamical behavior of the stochastic complex dynamical network. The numerical simulation of the GLV dynamics with empirical networks is shown in Fig. 2. We observe that the standard deviation of the realizations of the effective equation is always small (Fig. 2a), indicating that the drift term dominates the dynamics and the stochastic effect is not significant. Moreover, the standard deviation decreases after an initial jump (Fig. 2 c and e), indicating that the realizations of the stochastic complex dynamical network converge (Fig. 2 b and d). Therefore, the effective equation is a good approximation of the stochastic complex dynamical network.

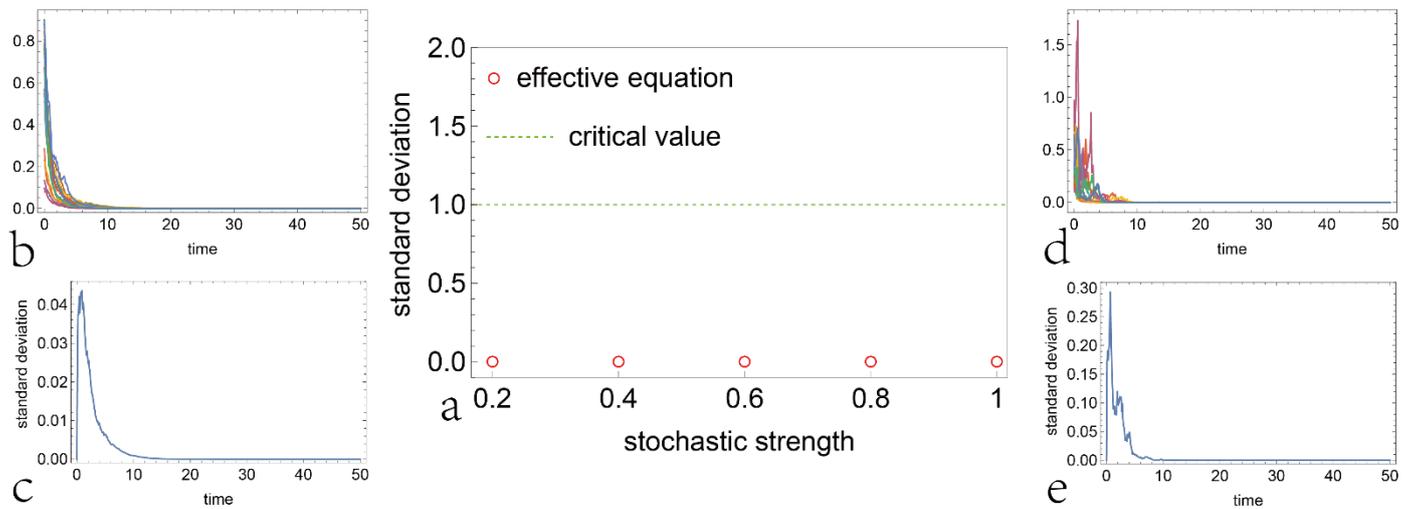

**Figure 2**: (a) The standard deviation of the realizations of the effective equation for different values of stochastic strength: 0.2, 0.4, 0.6, 0.8 and 1. (b-c) A sample realization of the stochastic complex dynamical network with stochastic strength 0.2 and the corresponding trajectory of the standard deviation of the effective equation. (d-e) A sample realization of the stochastic complex dynamical network with stochastic strength 1 and the corresponding trajectory of the standard deviation of the effective equation. We repeat each simulation 50 times.

# Discussion

Complex dynamical networks are ubiquitous in nature and engineering, but their analysis and control are hampered by their high dimensionality and the dependence of their dynamics on various factors, such as the topology of the interactions, the configuration of the parameters and the initial conditions of the nodes. Furthermore, most existing methods for dimensionality reduction of complex dynamical networks assume deterministic systems, while many real-world systems are subject to stochastic fluctuations. To overcome these challenges, we develop a general framework of dimensionality reduction for stochastic complex dynamical networks that can map a high-dimensional system with stochastic terms into a low-dimensional effective system with a single effective state variable and few effective parameters. Our framework is robust and versatile, as it can handle different types of dynamic models, network structures and node dynamics, and it can capture the collective dynamics of the stochastic complex dynamical networks, including the effects of determinism, stochasticity and convergence.

We have tested our framework of dimensionality reduction for stochastic complex dynamical networks on various examples, both synthetic and empirical, and have shown that it can accurately reproduce the collective dynamics of the original high-dimensional systems. In particular, we have demonstrated that the standard deviation of the realizations of the effective equation can capture the extent and trajectory of the collective dynamic behavior of the stochastic complex dynamical network. Our methodology is general and versatile, as it can be applied to a wide range of systems and contexts, such as biological, ecological, epidemiological and neural systems, and it can provide insights into the effects of determinism, stochasticity and convergence on the dynamics of complex networks.

# Data and code availability

The ready-to-use notebook codes to reproduce the results presented in the current study are available in OSF with the access code r7x3u (https://osf.io/r7x3u/).

# Acknowledgements

This work was supported by Zhejiang Provincial Natural Science Foundation of China (Grant No. LZ22G010001), Natural Science Fund of Zhejiang Province, China (LQ19G010004), and Science Foundation of Zhejiang Sci-Tech University (ZSTU) under Grant No. 18092125-Y and 22092034-Y.

# References

[1] Hughes TP, Kerry JT, Baird AH, Connolly SR, Dietzel A, Eakin CM, et al. Global warming transforms coral reef assemblages. Nature. 2018;556:492-6.

[2] Gauthier S, Bernier P, Kuuluvainen T, Shvidenko AZ, Schepaschenko DG. Boreal forest health and global change. Science. 2015;349:819-22.

[3] Huang H, Tu C, D'Odorico P. Ecosystem complexity enhances the resilience of plant-pollinator systems. One Earth. 2021;4:1286-96.

[4] Tu C, Suweis S, D'Odorico P. Impact of globalization on the resilience and sustainability of natural resources. Nature Sustainability. 2019.

[5] Schilders WHA, Van der Vorst HA, Rommes J. Model order reduction: theory, research aspects and applications: Springer; 2008.

[6] Brunton SL, Kutz JN. Data-driven science and engineering: Machine learning, dynamical systems, and control: Cambridge University Press; 2022.

[7] Jörg T. Generative Complexity in a Complex Generative World: A Generative Revolution in the Making: Springer Nature; 2021.

[8] Thurner S, Hanel R, Klimek P. Introduction to the theory of complex systems: Oxford University Press; 2018.

[9] Cheng X, Scherpen J. Model reduction methods for complex network systems. arXiv preprint arXiv:201202268. 2020.

[10] Lyapunov AM. The general problem of the stability of motion. International Journal of Control. 1992;55:531-4.

[11] Scheffer M, Bascompte J, Brock WA, Brovkin V, Carpenter SR, Dakos V, et al. Early-warning signals for critical transitions. Nature. 2009;461:53.

[12] Scheffer M, Carpenter SR, Lenton TM, Bascompte J, Brock W, Dakos V, et al. Anticipating critical transitions. Science. 2012;338:344-8.

[13] Suweis S, D'Odorico P. Early warning signs in social-ecological networks. PLoS ONE. 2014;9:e101851.

[14] Tu C, D'Odorico P, Suweis S. Critical slowing down associated with critical transition and risk of collapse in crypto-currency. Royal Society open science. 2020;7.

[15] Gao J, Barzel B, Barabási A-L. Universal resilience patterns in complex networks. Nature. 2016;530:307.

[16] Tu C, Grilli J, Schuessler F, Suweis S. Collapse of resilience patterns in generalized Lotka-Volterra dynamics and beyond. Phys Rev E. 2017;95:062307.

[17] Tu C, D'Odorico P, Suweis S. Dimensionality reduction of complex dynamical systems. iScience. 2021:101912.

[18] Laurence E, Doyon N, Dubé LJ, Desrosiers P. Spectral dimension reduction of complex dynamical networks. Physical Review X. 2019;9:011042.

[19] Vegué M, Thibeault V, Desrosiers P, Allard A. Dimension reduction of dynamics on modular and heterogeneous directed networks. arXiv preprint arXiv:220611230. 2022.

[20] Wu C, Duan D, Xiao R. A novel dimension reduction method with information entropy to evaluate network resilience. Physica A: Statistical Mechanics and its Applications. 2023:128727.

[21] Wilkinson DJ. Stochastic modelling for systems biology: CRC press; 2018.

[22] Black AJ, McKane AJ. Stochastic formulation of ecological models and their applications. Trends Ecol Evol. 2012;27:337-45.

[23] Chen W-Y, Bokka S. Stochastic modeling of nonlinear epidemiology. J Theor Biol. 2005;234:455-70.

[24] Laing C, Lord GJ. Stochastic methods in neuroscience: OUP Oxford; 2009.

[25] Battles Z, Trefethen LN. An extension of MATLAB to continuous functions and operators. SIAM Journal on Scientific Computing. 2004;25:1743-70.

[26] Trefethen LN. Approximation theory and approximation practice: Siam; 2013.

[27] Gardiner CW. Handbook of stochastic methods: springer Berlin; 1985.

[28] Holling CS. The components of predation as revealed by a study of small-mammal predation of the European pine sawfly. The Canadian Entomologist. 1959;91:293-320.

[29] Holling CS. The functional response of predators to prey density and its role in mimicry and population regulation. The Memoirs of the Entomological Society of Canada. 1965;97:5-60.

[30] Holling CS. Some characteristics of simple types of predation and parasitism. The Canadian Entomologist. 1959;91:385-98.

[31] Suweis S, Simini F, Banavar JR, Maritan A. Emergence of structural and dynamical properties of ecological mutualistic networks. Nature. 2013;500:449-52.

[32] Dakos V, Bascompte J. Critical slowing down as early warning for the onset of collapse in mutualistic communities. PNAS. 2014;111:17546-51.

[33] Suweis S, Grilli J, Banavar JR, Allesina S, Maritan A. Effect of localization on the stability of mutualistic ecological networks. Nat Commun. 2015;6:10179-.

[34] Cenci S, Song C, Saavedra S. Rethinking the importance of ecological networks through the glass of environmental variations. bioRxiv. 2017:219691.

[35] Grilli J, Adorisio M, Suweis S, Barabás G, Banavar JR, Allesina S, et al. Feasibility and coexistence of large ecological communities. Nat Commun. 2017;8.